\documentclass[11pt]{amsart}
\usepackage{amsmath, amsthm, amssymb, latexsym,url}
\usepackage{hyperref}

\author{J. Vandehey}
\thanks{Email: \href{mailto:vandehe2@illinois.edu}{\nolinkurl{vandehe2@illinois.edu}}}
\title[Normality in additive functions]{The normality of digits in almost constant additive functions}
\date{\today}
\keywords{Normal numbers, additive function, Selberg-Delange method}
\subjclass[2010]{Primary: 11K16}

\newtheorem{thm}{Theorem}[section]

\newtheorem{prop}[thm]{Proposition}

\begin{document}

\maketitle

\begin{abstract}
We consider numbers formed by concatenating some of the base $b$ digits from additive functions $f(n)$ that closely resemble the prime counting function $\Omega(n)$.  If we concatenate the last 
\[
\left\lceil y \frac{\log \log \log n }{\log b} \right\rceil
\] digits of each $f(n)$ in succession, then the number so created will be normal if and only if $0 < y \le 1/2$.  This provides insight into the randomness of digit patterns of additive function after the Erd\H{o}s-Kac theorem becomes ineffective.
\end{abstract}

\section{Introduction}\label{section:introduction}

We say a real number $z$ with fractional part $\{z \} = 0.z_1z_2 z_3\dots$ written in base $b$ is normal if every digit string $a_1 a_2 \dots a_k$ consisting of $k$ base $b$ digits occurs with relative frequency $b^{-k}$---in other words,
\[
\lim_{x\to \infty} 	\frac{\# \left\{1\le  n \le x-k \mid z_{n-1+i}=a_i, 1 \le i \le k\right\} }{x}= b^{-k}.
\]
(We shall, for the remainder of the paper, assume that the base $b$ is always the same unless otherwise specified, so we may save the reader from constant repetitions of the phrase ``in base $b$.'')

As a simple consequence of Birkhoff's pointwise ergodic theorem, almost all real numbers are normal \cite[Ch.~3]{DK}; however, despite their omnipresence among the reals, all numbers currently known to be normal have been explicitly constructed to be normal.  We do not know if any of the common mathematical constants such as $\pi$, $e$, or even $\sqrt{2}$ is normal.  

The first and simplest construction of a normal number is that of Champernowne \cite{champernowne}: in base $10$, the number
\[
0.(1)(2)(3)(4)(5)(6)(7)(8)(9)(10)(11)\dots
\]
composed by concatenating all natural numbers sequentially is normal.  Here we use the notation $(n)$ to refer to the string of digits that compose $n$, and the notation $(n)(m)$ to denote the concatenation of the two digit strings.

Not long thereafter, Copeland and Erd\H{o}s \cite{EC} proved that, in base $10$, the number
\[
0.(2)(3)(5)(7)(11)(13)(17)\dots
\]
composed by concatenating all primes is normal.  The proof of this fact is one of the few places in mathematics where the phrase, ``Because the primes are sufficiently dense in $\mathbb{N}$,'' is used in earnest.

From here, the study of normal numbers has flourished in a number of different settings.  Many mathematicians have sought to generalize the notion of normality to other settings and these have provided interesting results for continued fraction expansions \cite{AKS}, expansions with respect to the bases $-n \pm i$ \cite {M1}, matrix number systems \cite{M2}, $Q$-Cantor set expansions \cite{Mance1,Mance2,Mance3}, Markov shifts and intrinsically ergodic subshifts \cite{SM}, and higher-dimensional sequences (see \cite{LS} and the papers cited therein).

Studies into constructions for real normal numbers have tended towards one of two primary paths.  In one direction, several papers have investigated the base $b$ normality of numbers of the form
\[
\sum_{i=1}^\infty \frac{1}{b^{m_i}c^{n_i}}
\]
for special sequences $\{m_i\}$ and $\{n_i\}$.  The most powerful results of this type known to the author are those of Bailey and Crandall \cite{BC}; in their paper, they include a general conjecture that, if true, would prove the normality of a wide swath of common constants at once, including $\pi$.

We are more interested in the second direction, which follows in the footsteps of Champernowne, Copeland, and Erd\H{o}s.  These results are concerned with the normality of the numbers
\begin{align*}
\theta_f &= 0.(f(1))(f(2))(f(3))(f(4))(f(5))\dots\\
\tau_f &= 0.(f(2))(f(3))(f(5))(f(7))(f(11))\dots
\end{align*}
for some function $f(n)$.  Davenport and Erd\H{o}s \cite{DE} showed that $\theta_f$ is normal when $f(n)$ is a positive, integer-valued polynomial.  Nakai and Shiokawa \cite{NS} showed that both $\theta_f$ and $\tau_f$ are normal when $f(n)$ is the floor of a polynomial with real coefficients.  Madritsch, Thuswaldner, and Tichy \cite{MTT} likewise showed that both $\theta_f$ and $\tau_f$ are normal when $f(n)$ is the floor of an entire function with small logarithmic order.  De Konick and K\'{a}tai \cite{DKK1,DKK2} have investigated more number-theoretic functions $f(n)$ and have proved the normality of $\theta_f$ and $\tau_f$ for $f(n)=P(n+1)$, where $P(n)$ is the largest prime divisor of $n$.

In this paper, we want to extend these results to better understand the normality of $\theta_f$ when we let $f(n)$ be an additive function.  A function $f(n)$ is said to be additive if $f(mn)=f(m)+f(n)$ for relatively prime $m$ and $n$; such a function is defined entirely by its values on prime powers.  We will focus on additive functions closely related to the prime counting function $\Omega(n)$ and distinct prime counting function $\omega(n)$, which are defined by $\Omega(p^k)=k$ and $\omega(p^k)=1$.

However, it should come as no surprise that the constants $\theta_\Omega$ and $\theta_\omega$ are not normal.  The famous Erd\H{o}s-Kac theorem states that the values of $\omega(n)$ for $n\le x$ are normally distributed with mean and variance $\log \log x$ (and various generalizations state the same for $\Omega(n)$); thus, for almost all $n \le x$, given large $x$, the first digits of $\omega(n)$ should closely resemble the first half of the digits of $\log \log x$.  So if, say, $( \log \log x)$ contains an overabundance of zeros in the first half, we should expect the same to be true of $(\omega(n))$ for most $n\le x$.

The Erd\H{o}s-Kac theorem does not say anything about what happens in the last 
\[
\approx \frac{\log \log \log x}{2 \log b}
\]
digits of $\omega(n)$ (for $n\le x$).  To consider this further, we need a few extra definitions.  Let the truncation function $T_b(z,m)$ denote the string of the last $m$ base $b$ digits of $\lfloor z \rfloor$; we assume $z$ is always at least $0$.  If $\lfloor z\rfloor$ has fewer than $m$ digits, the truncation function $T_b(z,m)$ will include additional zeros at the head of the string to make sure the string returned is $m$ digits long.  Therefore,
\[
T_{10}(151,2) = (51) \qquad T_{10}(1,2)=(01) \qquad T_{10}(.5,2)= (00).
\]
Also, for $y>0$, let
\[
K_y(x) = \begin{cases} \left\lceil y\dfrac{\log \log \log x}{\log b} \right\rceil, & x > e^e,\\
1, & \text{otherwise}.
\end{cases}
\]
and let $\mathcal{K}(x)=K_{1/2}(x)$.   Alternately, the number $\mathcal{K}(x)$ may be roughly interpreted as half the digits of $\log \log x$ in base $b$, rounded up.

Now, given an additive function $f(n)$ and a fixed base $b$, let $f_y(n)= T_b(f(n),K_y)$ and let
\[
\theta_{f,y} = 0.(f_y(1))(f_y(2))(f_y(3))(f_y(4))(f_y(5))\dots.\footnote{Even though $f_y(n)$ is already considered to be just a string of digits, we include the extra parenthesis to prevent confusion.}
\] 
By the same argument we made above, Erd\H{o}s-Kac suggests that $\theta_{\omega,y}$ should not be normal for any $y>1/2$, but says nothing about any smaller $y$.  For that, we present the following new result.

\begin{thm}\label{thm:second}
For any $y$ satisfying $0<y \le 1/2$, the numbers $\theta_{\Omega,y}$ and $\theta_{\omega,y}$ are normal.  For any $y$ satisfying $1/2 < y$, the numbers $\theta_{\Omega,y}$ and $\theta_{\omega,y}$ are not normal. 
\end{thm}

So while the Erd\H{o}s-Kac theorem states roughly that the first half of the digits of $\omega(n)$ for $n\le x$ follow a very clear pattern, Theorem \ref{thm:second} states that the second half are statistically random.   

We wish to extend these results to a slightly larger class of functions.  We say a non-negative additive function $f$ is almost constant on primes if the following holds: there exist $c>0$, $1/2>\delta>0$, such that for all small $\epsilon>0$, we have 
\[
\exp\left( \sum_{p\le x} \frac{ B_\epsilon(x,p)}{p^{1-\delta}} \right) = o\left( \log x \right), \qquad x\to\infty,
\]
where
\[
B_\epsilon(x,p) = \min\left\{ 2, \dfrac{1}{(\log \log x)^{\epsilon/4}}|f(p)-c|\right\}
\]
For example, any additive function that satisfies $|f(p)-c|<p^{-\delta}$ for some $c,\delta>0$ is almost constant on primes.  

Likewise, we say that a non-negative additive function is weakly additive to mean the following.  For all small $\epsilon >0$, we have 
\[
\prod_{p\le x} \left( 1 - \sum_{p^k\le x, \ k\ge 2} \frac{C_\epsilon(x,p,k)}{p^k} \right) = 1+o(1), \qquad x\to \infty,
\]
where
\[
C_\epsilon(x,p,k) = \min\left\{ 2, \dfrac{1}{(\log \log x)^{1+\epsilon/4}} |f(p^k)-f(p^{k-1})-f(p)| \right\}.
\]
All completely additive functions, such as $\Omega(n)$, are weakly additive, as are all additive functions where $|f(p^k)-f(p^{k-1})-f(p)|$ always belongs to a bounded set, such as $\omega(n)$. 

With these definitions we can state a general theorem that contains Theorem \ref{thm:second} as a specific case.

\begin{thm}\label{thm:main}
Let $f(n)$ be a non-negative additive function that is almost constant on primes.  If $0<y \le 1/2$, then $\theta_{f,y}$ is normal.  If $1/2 < y$ and $f(n)$ is also weakly additive, then $\theta_{f,y}$ is not normal.
\end{thm}

Somewhat surprisingly, in the case $0 < y \le 1/2$, when $\theta_{f,y}$ is normal, we have no restrictions whatsoever on the values of $f$ at prime powers (outside of being non-negative).

The reason why $\theta_{f,y}$ fails to be normal is again due to the digits of $f(n)$ for $n\le x$ too closely resembling the digits of $\log \log x$, so one may wonder whether an Erd\H{o}s-Kac type distribution is the culprit.  Given an additive function, let 
\[
A_n = \sum_{p < n} \frac{f(p)}{p} \quad \text{ and } \quad B_n = \sum_{p < n } \frac{f(p)^2}{p}.
\]
There are various different conditions on $f(n)$ that guarantee it will be normally distributed about $A_n$ with variance $B_n$.  The simplest such conditions merely require that $B_n$ goes to infinity as $n$ goes to infinity and that $f(p)$ is uniformly bounded \cite{EK}.  Alternately, one may require that $B_n$ goes to infinity as $n$ goes to infinity and that
\[
\lim_{n\to \infty} B_n^{-1} \sum_{\substack{p < n\\ |f(p)| > \epsilon B_n^{1/2}}} \frac{f(p)^2}{p} = 0
\]
for every $\epsilon >0$ \cite{shapiro}.

However, there is nothing in the definition of $f(n)$ being almost constant on primes which implies either of these sets of conditions: on a very sparse set of primes, $f(p)$ could be enormous, and these few primes might be the primary contributers to $B_n$.  It is certainly possible that by ignoring such large values, we could show that being constant on primes would imply some type of Erd\H{o}s-Kac distribution on average.

For more on these types of results, see \cite{DH}.

\section{Bounds on exponential sums}\label{section:selbergdelange}

In the proof of normality results dealing with $\theta_f$ and $\tau_f$, a key technique is to reduce the problem to one concerning exponential sums, often of the form
\[
\sum_n e\left( \frac{\nu}{b^m} f(n) \right)
\]
for some integer $\nu$ and positive integer $m$, where the function $e(z)$ denotes $e^{2\pi i z}$.  When Davenport and Erd\H{o}s examined the case $f(n)$ is a polynomial, they used estimates on Weyl sums.  In our case, we want to examine the case when $f(n)$ is a non-negative additive function and will in turn use the Selberg-Delange method to estimate these sums.

We use results and terminology from Tenenbaum's book \cite[\S II.5.3]{tenenbaum}, mildly simplified to remove an unneeded parameter.  (In particular, we take $\delta=1$ in Tenenbaum's notation and replace his $c_0$ with our $\delta$.)

Let $z \in \mathbb{C}$, $\delta > 0$, $M >0$.  Then the Dirichlet series 
\[
F(s)=  \sum_{n=1}^\infty a_n n^{-s}
\]
is said to possess the property $\mathcal{P}(z;\delta,M)$ if the Dirichlet series 
\[
G(s;z) := F(s) \zeta(s)^{-z}
\]
can be extended to a holomorphic function on 
\[
\sigma \ge 1 - \delta/(1+\log^+ |\tau|), \quad (s=\sigma+i\tau)
\]
which is bounded on this domain in the following way 
\[
|G(s;z)| \le M.
\]

Suppose $F(s)$ has the property $\mathcal{P}(z;\delta,M)$ and there exist real positive numbers $\{b_n\}_{n=1}^\infty$ such that $|a_n| \le b_n$ for all $n$ and the series 
\[
\sum_{n=1}^\infty b_n n^{-s}
\]
satisfies the property $\mathcal{P}(w;\delta,M)$ for a complex number $w$, then $F(s)$ is said to be of the type $\mathcal{T}(z,w;\delta,M)$.

\begin{thm}\label{thm:selbergdelange}(\cite[\S II.5.3, Theorem 3]{tenenbaum} with $N=0$)
 Suppose the Dirichlet series $F(s) :=  \sum_{n=1}^\infty a_n n^{-s}$ is of the type $\mathcal{T}(z,w;\delta,M)$.  Then for $x \ge 3$, $A>0$, $|z| \le A$,  and $|w| \le A$, we have that 
\[
\sum_{n \le x} a_n = 
x(\log{x})^{z-1} \left( \frac{G(1;z)}{\Gamma(z)}+O\left(M\left(e^{-c_1 \sqrt{\log{x}}} + \frac{c_2}{\log{x}}\right)\right)\right)
\]
where the positive constants $c_1, c_2$ and the implicit constant in the big-O term depend on at most  $\delta$ and $A$.
\end{thm}

We now apply this result to the case we will need in our proof.

\begin{prop}\label{prop:main}
Let $f(n)$ be a non-negative additive function that is almost constant on primes, with constant $c$.  Let $\eta$ and $\epsilon$ be fixed positive real numbers strictly less than $1$, and let $a$ be a non-zero integer satisfying $|a| < \eta^{-2}$.
\begin{enumerate}
\item For $\epsilon \mathcal{K}(x) \le m \le (1-\epsilon) \mathcal{K}(x)$, we have
\[
\sum_{n\le x} e\left(\frac{a}{b^m} f(n)\right) = o(x), \quad x\to \infty,
\]
where the rate of decay in the little-o term is uniform for fixed $\eta$ and $\epsilon$.  

\item If $f(n)$ is also weakly additive, then for $ (1+\epsilon) \mathcal{K}(x)\le m$, we have
\[
\sum_{n\le x} e\left(\frac{a}{b^m} f(n)\right)  = x e\left(\frac{a}{b^m} c \log \log x\right) (1+o(1)), \quad x \to \infty
\]
where the rate of decay in the little-o term is uniform for fixed $\eta$ and $\epsilon$.  
\end{enumerate}
\end{prop}

Recall, by the definition of $K_y(x)$ that
\[
\mathcal{K}(x) := \begin{cases}
\left\lceil \dfrac{\log \log \log x}{2 \log b} \right\rceil, & x> e^e, \\
1, & \text{otherwise}.
\end{cases}
\]

\begin{proof}
For a fixed $x$, we shall consider the function $F(s)$ with coefficients $a_n$, which are assumed to be multiplicative and defined on prime powers by
\[
a_{p^k}:= \begin{cases} e\left(\dfrac{a}{b^m} f(p^k) \right) & p^k\le x\\
e\left(\dfrac{a}{b^m} (f(p^{k-1})+f(p)) \right) & p \le x < p^k\\
e\left(\dfrac{a}{b^m} k c \right)& p > x\end{cases}.
\]

Let $c'=e(ac/b^m)$.  Then to the right of the line $\sigma=1$, we have
\begin{align*}
G(s;c')&= F(s) \zeta(s)^{-c'}\\
&=\prod_p \left( 1+ \frac{a_p}{p^s}+\frac{a_{p^2}}{(p^2)^s} + \dots \right)  \left( 1- \frac{1}{p^s} \right)^{c'}\\
&=\exp\left( \sum_{p} \log  \left( 1+ \frac{a_p}{p^s}+\frac{a_{p^2}}{(p^2)^s} + \dots \right)  + c' \log   \left( 1- \frac{1}{p^s} \right) \right)\\
&= \exp \left( \sum_{p \le x } \frac{a_p-c'}{p^s}  + \sum_p O(p^{-2s}) \right).
\end{align*}
This shows that $G(s,c')$ can be extended to a holomorphic function on $\sigma > 1/2$.  

The difference $a_p-c'$ can be bounded in one of two ways: first, by
\begin{align*}
|a_p-c'|&=\left| e\left(\frac{ac}{b^m}\right)\left( e\left(\frac{a(f(p)-c)}{b^m}\right)-1 \right) \right| \\
&\le \frac{2\pi a |f(p)-c|}{b^m}\\
& \le \frac{ 2\pi \eta^{-2} |f(p)-c|}{(\log \log x)^{\epsilon/2}},
\end{align*}
using $|e^{i x} -1|\le |x|$, and, second, by $2$ since $|a_p|=|c'|=1$.

Provided $x$ is sufficiently large so that $(\log \log x)^{\epsilon/4}\ge 2 \pi\eta^{-2}$, we have that 
\[
 \exp \left( \sum_{p \le x } \frac{a_p-c'}{p^s} \right) = o(\log x)
\]
for $\sigma \ge 1-\delta $ by the definition of $f$ being almost constant on primes.  Thus, for some $M=M(x)$ that is $o(\log x)$, we have that $G(s;c')\le M$ for $\sigma \ge 1-\delta$, and so $F(x)$ has property $\mathcal{P}(c';\delta,M)$.

Since each $|a_n|=1$, we can let $b_n=1$ and note that the zeta function itself has property $\mathcal{P}(1;\delta,M)$ rather trivially, so that $F(s)$ has property $\mathcal{T}(c',1;\delta, M)$.  We can therefore apply Theorem \ref{thm:selbergdelange}.

Thus, we have that
\[
\sum_{n\le x} e(af(n)/b^m) = x(\log x)^{c'-1} \left( \frac{G(1;c')}{\Gamma(c')} + o(1)\right)
\]
by the bound on $M$.  The function bounded by $o(1)$ here decays to zero as $x$ tends to infinity and depends only on $\epsilon$ and $\eta$.

Suppose we are in the case where $\epsilon \mathcal{K} (x)\le m \le (1-\epsilon) \mathcal{K}(x)$.  Then the real part of $c'-1$ is bounded by
\begin{align*}
\operatorname{Re}(e(ac/b^m)-1)&\le -\frac{1}{2!}\left(\frac{2\pi ac}{b^m}\right)^2+ \frac{1}{4!}\left(\frac{2\pi ac}{b^m}\right)^4\\
&\le -\frac{1}{2}\left(\frac{2\pi ac}{b^{(1-\epsilon)\mathcal{K}}}\right)^2(1+O(( \log \log x)^{-\epsilon}))\\
&\le -e^{O(1)} (\log\log x)^{-(1-\epsilon)}(1+O(( \log \log x)^{-\epsilon})).
\end{align*}
Therefore, we have
\[
x(\log x)^{c'-1} \le x\exp\left(-e^{O(1)} (\log\log x)^{\epsilon}(1+O(( \log \log x)^{-\epsilon})) \right) = o(x).
\]

The function $1/\Gamma(z)$ is uniformly bounded in any sufficently small compact neighborhood of $1$, and $c'$ is in a small neighborhood of $1$ for suffficiently large $x$.

To finish the proof of the first half of the proposition, it suffices to show that $G(1;c')=O(1)$, which itself would be implied by
\[
\sum_{p\le x} \frac{B_\epsilon(x,p)}{p} = O(1)
\]
for all sufficiently small $\epsilon$.  By the definition of $f(n)$ being almost constant on primes, we know that
\[
\sum_{p\le x} \frac{B_\epsilon(x,p))}{p^{1-\delta}} = O(\log \log x),
\] 
and we also have that $B_\epsilon(x,p)$ are positive, decreasing functions in $x$.  Therefore, by partial summation, we have
\begin{align*}
\sum_{p\le x} \frac{B_\epsilon(x,p)}{p} &= x^{-\delta} \sum_{p\le x} \frac{B_\epsilon(x,p)}{p^{1-\delta}} +\delta \int_1^x t^{-\delta-1} \sum_{p\le t} \frac{B_\epsilon(x,p)}{p^{1-\delta}} \ dt\\
&\le O(x^{-\delta} \log \log x) + \delta \int_1^x t^{-\delta-1} \sum_{p\le t} \frac{B_\epsilon(t,p)}{p^{1-\delta}} \ dt \\
&= O(x^{-\delta} \log \log x) + O\left( \int_1^x t^{-\delta-1} \log \log t \ dt \right)\\
&=O(1),
\end{align*}
which completes the proof of statement ($1$).

Now suppose we are in the case where $(1+\epsilon) \mathcal{K}(x)\le  m$ and $f(n)$ is weakly additive.  Then, in this case, we have
\[
c'-1 =2\pi i \frac{a }{b^m} c+ O\left(\frac{1}{b^{2m}}\right),
\]
so that
\begin{align*}
(\log x)^{c'-1} &= e\left(\frac{a}{b^m} c \log \log x\right) \cdot \exp\left(O\left( \frac{\log \log x}{b^{2m}} \right) \right)\\
&= e\left(\frac{a}{b^m} c \log \log x\right)  (1+O((\log \log x)^{-\epsilon})).
\end{align*}
Also, since $c'-1=O((\log \log x)^{-1})$, we have
\[
\frac{1}{\Gamma(c')} = 1+O((c'-1)) = 1+O((\log \log x)^{-1}).
\]

Therefore, to finish the proof in this case, it suffices to show that $G(1;c')=1+o(1)$ uniformly in $x$.  Let $G(z)$ denote the function
\[
\prod_p \left( 1-\frac{z}{p}\right)^{-1} \left( 1-\frac{1}{p}\right)^z.
\]
We have $G(1)=1$ and $G(z)$ is analytic at $1$, so that $G(c')=1+O((\log \log x)^{-1})$.  Moreover, we can factor such a term out of $G(1;c')$ to leave
\[
G(1;c') = G(c') \prod_{p\le x} \left( 1+ \sum_{k=1}^\infty \frac{a_{p^k}}{p^k} \right) \left( 1- \frac{a_p}{p} \right).
\]
(We do not need that the product is absolutely convergent since we are at most changing finitely many factors.)  Now we examine this latter product in greater detail.  We have, first, that
\begin{align*}
 \prod_{p\le x} \left( 1+ \sum_{k=1}^\infty \frac{a_{p^k}}{p^k} \right) \left( 1- \frac{a_p}{p} \right)&= \prod_{p\le x} \left( 1+ \sum_{k \ge 2} \frac{1}{p^k}(a_{p^k}-a_{p^{k-1}}a_p) \right) \\
&=  \prod_{p\le x} \left( 1+ \sum_{p^k \le x, \ k \ge 2} \frac{1}{p^k}(a_{p^k}-a_{p^{k-1}}a_p) \right),
\end{align*}
and provided $x$ is sufficiently large, we have also, by a similar argument to earlier in the proof, that this is bounded above and below by 
\[
\prod_{p\le x}  \left( 1+ \sum_{p^k \le x, \ k \ge 2} \frac{C_\epsilon(x,p,k)}{p^k} \right) \quad \text{and} \quad  \prod_{p\le x}  \left( 1- \sum_{p^k \le x, \ k \ge 2} \frac{C_\epsilon(x,p,k)}{p^k} \right)
\]
respectively.  The latter equals $1+o(1)$, since $f(n)$ is weakly additive, and this in turn implies that the former is $1+o(1)$ as well, since the logarithm of each factor in the latter product is larger in norm than the logarithm of the corresponding factor in the former product.  Thus we have $G(1;c')=1+o(1)$ and this completes the proof.
\end{proof}

If we only care about the case $\epsilon \mathcal{K}(x) \le m \le (1-\epsilon) \mathcal{K}(x)$, then the above proof still works under weaker conditions on $f(n)$.  For example, we could replace $o(\log x)$ in the definition of $f(n)$ being almost constant on primes with any function that would satisfy $o((\log x)^{2-c'})$, such as $O(\log x \log \log x)$.

Similarly, if we consider an alternate definition of the truncation function $T^*_\epsilon(z,x)$ which returns the string of the last $(1-\epsilon)\mathcal{K}(x)$ through $\epsilon\mathcal{K}(x)$ digits of $\lfloor z \rfloor$, and let $f^*_\epsilon(n)=T^*_\epsilon(f(n),n)$, then we would expect the number
\[
\theta_{f,\epsilon}^* = 0.(f^*_\epsilon(1))(f^*_\epsilon(2))(f^*_\epsilon(3))(f^*_\epsilon(4))\dots
\]
to be normal while only requiring the definition of being almost constant on primes to hold for this specific $\epsilon$.  This suggests that the ``easiest'' place for randomness to occur in the digits of an additive function is around the 
\[
\frac{1}{4} \frac{\log \log \log x}{\log b} \text{th}
\]
place.

\section{Proof of Theorem \ref{thm:main}}

Consider a string of $k$ base $b$ digits $a_1 a_2 \cdots a_k$.  Let $N^*(x)$ denote the number of times the string occurs in the first $x$ digits of $\theta_{f,y}$ after the decimal point.  To prove Theorem \ref{thm:main}, it suffices to prove
\begin{equation*}
\lim_{x\to \infty} \frac{1}{x} N^*(x)\ \begin{cases}
=b^{-k}, & \text{if }0 < y \le 1/2, \\
\neq b^{-k}, & \text{if }1/2 < y.
\end{cases}
\end{equation*}

The $x$th digit of $\theta_{f,y}$ may occur in the middle of some string $(f_y(x'))$.  If we let $N(x)$ denote the number of times we see the string occur within the digits of 
\[
(f_y(1))(f_y(2))(f_y(3))\dots (f_y(x')),
\] 
then it is clear that $N^*(x)=N(x')+O(K_y(x'))$, since the latter function adds at most $K_y(x')$ new places for the string to occur.  Moreover, the string $(f_y(n))$ is exactly $K_y(n)$ digits long, so, by a simple application of partial summation, the entire string 
\[
 (f_y(1))(f_y(2))(f_y(3))\dots (f_y(x')),
\]
is $x' K_y(x')(1+o(1))$ digits long.  Therefore, we have
\[
\frac{1}{x} N^*(x) = \frac{1}{x'K_y(x')} N(x')(1+o(1))
\]
and so it suffices now to prove
\begin{align*}
\lim_{x\to \infty} \frac{1}{xK_y(x)} N(x)&=b^{-k}, \qquad \text{if }0 < y \le 1/2, \\
\limsup_{x\to \infty} \frac{1}{xK_y(x)} N(x) &> b^{-k}, \qquad \text{if }1/2 < y.
\end{align*}

We can rewrite $N(x)$ in the following way:
\begin{equation}\label{eq:Nx}
N(x) = \sum_{n=1}^x \sum_{m=k}^{K_y(n)} \theta(b^{-m} f(n)) + O(x)
\end{equation}
where
\[
\theta(z):=\begin{cases}
                                        1, & 0.a_1 a_2 \cdots a_k \le \{z\} < 0.a_1 a_2 \cdots a_k + b^{-k}, \\
                                        0, & \text{otherwise}.
                                      \end{cases}
\]
The term $\theta(b^{-m} f(n))$ is an indicator of whether the string $a_1a_2 \dots a_k$ occurs in $f(n)$ starting at the $m$th place (to the left of the decimal).  The big-O term accounts for occurences of the string that start in a given $f_y(n)$ and finish in another.

As noted above, the whole string
\[
 (f_y(1))(f_y(2))(f_y(3))\dots (f_y(x')),
\]
contains $xK_y(x)(1+o(1))$ digits, so in \eqref{eq:Nx} we may replace $K_y(n)$ with $K_y(x)$ with an error of $o(xK_y(x))$ and switch the order of summation to obtain 
\begin{equation}\label{eq:Nx2}
N(x) = \sum_{m=k}^{K_y(x)} \sum_{n=1}^x  \theta(b^{-m} f(n)) + o(xK_y(x)).
\end{equation}

Let $Y=\min\{y,1/2\}$.  Then for any small, fixed $\epsilon >0$, let
\begin{align*}
U(x) &=U(x,\epsilon) = \sum_{m=\epsilon \mathcal{K}(x)}^{(1-\epsilon)K_Y(x)} \sum_{n=1}^x  \theta(b^{-m} f(n)) \\
V(x) &=V(x,\epsilon) = \sum_{m=(1+\epsilon)K_Y(x)}^{K_y(x)} \sum_{n=1}^x  \theta(b^{-m} f(n)) .
\end{align*}
Note that 
\[
\mathcal{K}(x)=\frac{1/2}{y}K_y(x)+O_y(1),
\]
and that unless $1/2 <y$, the sum $V(x)$ will be identically zero.

Thus we may write $N(x)$ as 
\[
N(x) =U(x)+V(x) + O(\epsilon x K_y(x)) + o(xK_y(x)).
\]

We estimate $U(x)$ first.  Following Davenport and Erd\H{o}s, we pick a small positive constant $\eta$ and consider two functions $\theta_1(z)$ and $\theta_2(z)$ both periodic in $z$ with period 1, with $\theta_1(z) \le \theta(z) \le \theta_2(z)$ and having Fourier expansions
\begin{align*}
\theta_1(z) &= (b^{-k}- \eta) + \sum_{\nu \neq 0} A_\nu^{(1)} e(vz) \\
\theta_2(z) &= (b^{-k}+ \eta) + \sum_{\nu \neq 0} A_\nu^{(2)}e(vz) \end{align*}
where the coefficients are bounded by 
\[
|A_\nu^{(*)}| \le \min \left( \frac{1}{|\nu|},\frac{1}{\eta\nu^2}\right).
\]
Note that this makes the sums absolutely convergent and hence we can change the order of summation without trouble.  (For the existence of such functions, see \cite[pp.~91--92,99]{koksma}.)

In particular, this gives, for $Y=\min\{y,1/2\}$ as above,
\begin{align*}
U(x) &\ge (b^{-k}-\eta) x((1-\epsilon)K_Y(x)-\epsilon\mathcal{K}(x)) +  \sum_{m=\epsilon \mathcal{K}(x)}^{(1-\epsilon)K_Y(x)} \sum_{n=1}^x \sum_{\nu \neq 0} A_\nu^{(1)} e\left( \frac{\nu }{b^m}f(n)\right)\\
&=b^{-k}xK_Y(x)(1+O(\epsilon)+O(\eta)) +   \sum_{\nu \neq 0} A_\nu^{(1)} \sum_{m=\epsilon \mathcal{K}(x)}^{(1-\epsilon)K_Y(x)} \sum_{n=1}^x  e\left( \frac{\nu }{b^m} f(n)\right)\\
&= b^{-k}xK_Y(x)(1+O(\epsilon)+O(\eta))  +O(x K_y(x) \eta) \\
&\qquad + \sum_{\substack{\nu \neq 0\\ |\nu| \le \eta^{-2}}} A_\nu^{(1)} \sum_{m=\epsilon \mathcal{K}(x)}^{(1-\epsilon)K_Y(x)} \sum_{n=1}^x  e\left( \frac{\nu }{b^m} f(n)\right).
\end{align*}
Applying Proposition \ref{prop:main}, we then obtain
\begin{align*}
U(x)&\ge b^{-k}xK_Y(x)(1+O(\epsilon)+O(\eta)) + \sum_{\substack{\nu \neq 0\\ |\nu| \le \eta^{-2}}} A_\nu^{(1)} \sum_{m=\epsilon \mathcal{K}(x)}^{(1-\epsilon)K_Y(x)} o(x)\\
&= b^{-k} xK_Y(x)(1+O(\epsilon)+O(\eta)+o_{\eta,\epsilon}(1)),
\end{align*}
where the little-o term tends to zero uniformly for fixed $\eta$ and $\epsilon$.
By applying the same idea with $\theta_2(z)$, we obtain
\[
U(x) \le  b^{-k} xK_Y(x)(1+O(\epsilon)+O(\eta)+o_{\eta,\epsilon}(1))
\]
as well.

We can now complete the proof in the case $0 < y \le 1/2$.  In this case, we have
\[
\frac{1}{xK_y(x)} N(x) = b^{-k}(1+O(\epsilon)+O(\eta)+o_\eta(1)).
\]
We have chosen $\epsilon$ and $\eta$ to be fixed as $x$ varies, so this does not immediately give us that the limit is $b^{-k}$ (because we have no guarantee that the limit exists).  However, it does show that the lim sup and lim inf of $N(x)/xK_y(x)$ are both equal to $b^{-k}(1+O(\epsilon)+O(\eta))$; by taking $\epsilon$ and $\eta$ arbitrarily small, we see that the lim sup is at most $b^{-k}$ and the lim inf is at least $b^{-k}$, so the full limit itself must exist and equal $b^{-k}$.

Now we investigate the case $1/2 < y$.  Consider a function $\theta_3(z)$ that satisfies the following properties:
\begin{enumerate}
\item $0 \le \theta_3(z) \le \theta(z)$;
\item $\theta_3(z)=\theta(z)$ except on the intervals 
\[
[0.a_1a_2 \dots a_k-b^{-2 K_y(x)},0.a_1a_2 \dots a_k+b^{-2K_y(x)}]
\]
 and 
\[
[0.a_1a_2 \dots a_k + b^{-k}-b^{-2 K_y(x)},0.a_1a_2 \dots a_k + b^{-k}+b^{-2K_y(x)}];
\] 
and,
\item $\theta_3(z)$ is continuous and piecewise smooth.
\end{enumerate}
In particular, since $\theta_3(z)$ is assumed to be continuous and piecewise smooth, its Fourier series converges absolutely \cite[p.~81]{tolstov}, so the order of summation doesn't matter.  We shall again refer to the Fourier coefficients by $A_\nu^{(3)}$.

Therefore, given some small $\eta$, we have
\begin{align*}
V(x) &\ge \sum_{m=(1+\epsilon) \mathcal{K}(x)}^{K_y(x)} \sum_{n=1}^x \sum_{\nu} A_\nu^{(3)} e\left( \frac{\nu }{b^m} f(n)\right)\\
&= \sum_{\nu} A_\nu^{(3)} \sum_{m=(1+\epsilon) \mathcal{K}(x)}^{K_y(x)} \sum_{n=1}^x e\left( \frac{\nu }{b^m} f(n)\right)\\
&= \sum_{|\nu|< \eta^{-2}} A_\nu^{(3)} \sum_{m=(1+\epsilon) \mathcal{K}(x)}^{K_y(x)} \sum_{n=1}^x e\left( \frac{\nu }{b^m}f(n) \right)+ O(\eta x(K_y(x)-(1+\epsilon)\mathcal{K}(x))).
\end{align*}
Applying Proposition \ref{prop:main} again and noting that we can add the $|\nu|\ge \eta^{-2}$ terms back in at the cost of another copy of
\[
O(\eta x(K_y(x)-(1+\epsilon)\mathcal{K}(x))),
\]
this becomes 
\begin{align*}
V(x)&\ge \sum_{|\nu|< \eta^{-2}} A_\nu^{(3)} \sum_{m=(1+\epsilon) \mathcal{K}(x)}^{K_y(x)} xe\left( \frac{\nu}{b^m} \log \log x\right)(1+o(1))\\
&\qquad + O(\eta x(K_y(x)-(1+\epsilon)\mathcal{K}(x)))\\
&= \sum_\nu A_\nu^{(3)} \sum_{m=(1+\epsilon) \mathcal{K}(x)}^{K_y(x)} xe\left( \frac{\nu}{b^m} \log \log x\right) \\
&\qquad +x (K_y(x)-(1+\epsilon)\mathcal{K}(x))(o(1)+O(\eta)).
\end{align*}
The sum
\begin{align*}
\sum_\nu A_\nu^{(3)} \sum_{m=(1+\epsilon) \mathcal{K}(x)}^{K_y(x)} e\left( \frac{\nu}{b^m} \log \log x\right)&=\sum_{m=(1+\epsilon) \mathcal{K}(x)}^{K_y(x)}  \sum_\nu A_\nu^{(3)}  e\left( \frac{\nu}{b^m} \log \log x\right)\\
&= \sum_{m=(1+\epsilon) \mathcal{K}(x)}^{K_y(x)} \theta_3\left( \frac{1}{b^m} \log \log x\right)
\end{align*}
approximates the number of times the string occurs in the digits of $\log \log x$ to the left of the $(1+\epsilon)\mathcal{K}$th place.  It will, in fact, equal and not just approximate the number of occurences of the string if we make two assumptions: that $x$ is sufficiently large and that $\log \log x$ is an integer ending in a non-zero digit base $b$.  In order for $\theta_3$ evaluated at $\log \log x / b^m$ to differ from $\theta$ at the same point, the fractional part of $\log\log x /b^m $ must be within $b^{-2K_y(x)}$ of a number in $[0,1]$ that has all zeros after the $k$th place, but since we have assumed that $\log \log x$ is an integer ending in a non-zero digit and since $m$ is at least $\mathcal{K}(x)$ but at most $K_y(x)$, the fractional part of $\log \log x / b^m$ has a non-zero digit after the $k$th place, but none after the $K_y(x)+1$th place.

Thus, for $x$ sufficiently large and $ \log \log x $ an integer ending on a non-zero digit, we have that 
\[
V(x) \ge x \sum_{m=(1+\epsilon) \mathcal{K}(x)}^{K_y(x)} \theta\left( \frac{1}{b^m} \log \log x\right) +x (K_y(x)-(1+\epsilon)\mathcal{K}(x))(o(1)+O_\eta(1)).
\]

Now let $M$ be a large integer and consider $x_M$ so that $\log \log x_M$ is an integer whose digit string is composed of $M$ copies of the string $a_1 a_2 \dots a_k$ concatenated together followed by exactly $\mathcal{K}(x)$ digits, the last of which is non-zero.  In particular, this implies that $M = (K_y(x_M)-\mathcal{K}(x_M))/k$.  In this case, we have
\begin{align*}
N(x_M) &\ge b^{-k}x_M\mathcal{K}(x_M) (1+O(\epsilon)+O(\eta)+o_\eta(1)) + x_M\frac{K_y(x_M)-\mathcal{K}(x_M)}{k} \\
&\qquad +x_M K_y(x_M)(O(\epsilon )+o(1)+O_\eta(1))\\
&= x_M K_y(x_M)\left( \frac{1}{k}\left(1-\frac{1}{2y}\right)+ \frac{1}{b^k}\frac{1}{2y}\right) (1+O(\epsilon)+O(\eta)+o_\eta(1)+O_\eta(1)).
\end{align*}
By choosing $\epsilon$ and $\eta$ sufficiently small, we see that the lim sup of $N(x)/xK_y(x)$ must exceed $b^{-k}$.

Thus, the limit of $N(x)/xK_y(x)$ cannot be $b^{-k}$, so $\theta_{f,y}$ is not normal in this case.

\section{Acknowledgements}

The author wishes to thank Heini Halberstam for his helpful comments.

\end{document}